\documentclass{amsart}
\usepackage{amssymb,amscd, hyperref, bm, color }
\usepackage{amsmath}
\usepackage{enumerate}
\usepackage{stmaryrd}  
\usepackage[mathcal]{eucal}
\usepackage{array,float}
\usepackage{longtable}
\usepackage{multirow}
\usepackage{cite}
\usepackage[left=3cm, right=3cm]{geometry}

\usepackage{xy}
\input xy
\xyoption{all}
\usepackage{pdflscape} 
\usepackage{hyperref}
\hypersetup{
	colorlinks=blue,
	linkcolor=blue,
	filecolor=magenta,      
	urlcolor=cyan,
}

\usepackage[draft]{graphicx}

\numberwithin{equation}{section}
\newtheorem{thm}{Theorem}[section]

\newtheorem{corollary}[thm]{Corollary}
\newtheorem{lemma}[thm]{Lemma}

\theoremstyle{definition}

\newtheorem*{remark*}{Remark}
\newtheorem{remark}[thm]{Remark}

\newcommand{\compcent}[1]{\vcenter{\hbox{$#1\circ$}}}

\newcommand{\comp}{\mathbin{\mathchoice
		{\compcent\scriptstyle}{\compcent\scriptstyle}
		{\compcent\scriptscriptstyle}{\compcent\scriptscriptstyle}}}

\newcommand{\SL}{\mathrm{SL}}
\newcommand{\GL}{\mathrm{GL}}

\newcommand{\mat}[4]{\left[ \begin{matrix}
		#1 & #2 \\  #3 & #4 \end{matrix}\right] }

\renewcommand{\det}{{\mathrm{det}}}

\newcommand{\Ind}{{\mathrm{Ind}}}

\title[Tensor product of irreducible characters of $\GL_2(\mathbb F_q)$]{Tensor product of irreducible characters of $\GL_2(\mathbb F_q)$}

\author{Archita Gupta}
\address{ AG:  Department of Mathematics and Statistics,
	Indian Institute of Technology Kanpur,
	Kanpur 208016, India. }
\email{architagup20@iitk.ac.in}

\author{M Hassain}
\address{MH: Harish-Chandra Research Institute, a CI of Homi Bhabha National Institute, Chhatnag Road, Jhunsi,
	Prayagraj - 211019, India.}

\email{hassainm@hri.res.in}



\begin{document}
	%\title{Tensor product of irreducible characters of $\GL_2(\mathbb F_q)$} 

	\begin{abstract}
		%We decompose the tensor product of two irreducible characters of $\GL_2(\mathbb F_q)$ for odd $q$ and classify pairs such that their tensor product is multiplicity free. We also show that  in each case is ... We also show that .... 
	We decompose the tensor product of two irreducible representations of $\GL_2(\mathbb F_q)$ for odd $q$ and classify the pairs such that their tensor product is multiplicity free. We also classify the pairs such that their tensor product has unique decomposition property. We additionally characterize the self-dual irreducible representations of $\GL_2(\mathbb F_q).$
	\end{abstract}
\maketitle
\vspace{-20pt}
\section{Introduction}

	Let $G$ be a finite group and $V_1$ and $V_2$ be any complex irreducible representations of $G.$ The problem of decomposing the tensor product $V_1 \otimes V_2$ into the irreducible constituents is called the "tensor product problem of the group $G$". This problem is known to be very hard in general. C. Bessenrodt and A. Kleschev \cite{MR1722888} proved that a tensor product of any two higher dimensional complex irreducible representations of $S_n$ is always reducible. In fact, they proved that such a tensor product is inhomogeneous; that is, the decomposition contains at least two distinct irreducible constituents. In the same article, they also classified the pairs $(V_1, V_2)$ of representations of $A_n$ for which $V_1 \otimes V_2$ is irreducible. The parallel questions for the modular representations of $S_n$ and $A_n$ were considered in \cite{MR1750169} and \cite{MR1764578}, and they classified the pairs $(V_1, V_2)$ such that $V_1 \otimes V_2$ has atmost three non-isomorphic irreducible constituents and conjectured for four constituents in \cite{MR1722888}.
	C. Bessenrodt and C.Bowman \cite{MR3720803} classified the pairs $(V_1, V_2)$ of complex irreducible representations of $S_n$ such that $V_1 \otimes V_2$ is multiplicity free.
	E. Letellier \cite{MR3022764} proved that given  unipotent characters $U_1,U_2,...,U_k$ of  $\mathrm{GL}_n(\mathbb{F}_q),$ $\langle U_1\otimes U_2\otimes  \cdots \otimes U_k , 1\rangle$ is  a polynomial in $q$ with non-negative integer coefficients. T. Hausel et al \cite{MR3003926},\cite{MR3034296} proved identical results for generic complex split semisimple characters of $\mathrm{GL}_n(\mathbb{F}_q).$ This study is generalised by removing the genericity constraint by T. Scognamiglio \cite{scognamiglio2023generalization}.
	%	
	%	For $G$, a finite simple group of Lie type, Gerard Heidi et al. have given a complete description of the constituents of the tensor product $\mathrm{St}\otimes \mathrm{St}$ of the Steinberg character $\mathrm{St}$ of $G$ in \cite{MR3056296}.
	%	
	%	In \cite{MR3003926,MR3034296}  authors study the multiplicities of trivial representation in the tensor products of generic complex split semisimple characters of $\mathrm{GL}_n(\mathbb{F}_q).$ This study is generalized by removing genericity constraint in \cite{scognamiglio2023generalization}.
	%
	%	For the general case, in \cite{MR4201484}, Gabriel Navarro and Pham hu tiep classified all finite quasisimple groups such that $V_1\otimes V_2$ is irreducible when $V_1$ and $V_2$ are of same degree.
	%	
	%	

	In this article, we consider the tensor product problem of $G=\GL_2(\mathbb F_q),$ where $q=p^n$ for an odd prime $p.$ Our main objective is to comprehend the complete decomposition of the tensor product of any two irreducible characters of $G.$ This problem is also considered in an unpublished work \cite{gurjyot}. The conclusions reached there are highly computational in nature because they rely on the character table of $G$ and require additional computation for the complete results. However, we proceed by taking a more conceptual approach and using the available ideas for the tensor product problem of $G.$ L. Aburto Hageman and Jose Pantoja \cite{MR1757476} reduced the decomposition of tensor products of $G$ as a direct sum of certain induced representations from the tori (see Theorem~\ref{thm:Jose Pantoja tensor product}). We mainly proceed by studying the irreducible constituents of these induced representations in Section~\ref{sec:decomp. of induction from tori}.
	
		Any irreducible representation of $G$ has dimension in $\{1,q-1,q,q+1\}$ (see for example \cite{MR1153249}, Chapter 5). We follow the notations of \cite{MR1757476} for the representations of $G.$ Let $\chi_\alpha^1, \chi_\alpha^q, \chi_{(\beta,\gamma)}^{q+1}$ and $\chi_{\Lambda}^{-1}$ denote one-dimensional, $q$-dimensional, $(q+1)$-dimensional and $(q-1)$-dimensional (cuspidal) irreducible representations of $G$ respectively, where $\alpha, \beta ,\gamma$ are characters of ${\mathbb{F}_q^\times}$ and $\Lambda$ is an indecomposable character of the anisotropic torus $T_{-1}$ of $G.$

	%For more details of the notations, see section~\ref{sec:decomp. of induction from tori}.
	
	% Let $\mathbb{F}_{q^2}^\times$ be the unique quadratic extension of $\mathbb{F}_q^\times.$ Let $\bar{\Lambda}$ denote the restriction of the character $\Lambda$ of $\mathbb{F}_{q^2}^\times$ to $\mathbb{F}_q^\times.$ 

	For $ \alpha \in \widehat{\mathbb{F}_q^\times},$
	we  define the following sets containing inequivalent irreducible representations of $G$:
	\begin{eqnarray*}
		\mathbb{S}_\alpha &=& \{  \chi^{q+1}_{(\gamma_1, \gamma_2)}  \mid \gamma_1\neq \gamma_2 , \alpha=\gamma_1\gamma_2 \},\\
		\mathbb{W}_\alpha &=&  \{ \chi_\Lambda^{-1} \mid \bar{\Lambda}=\alpha, 
		\Lambda \,\, \mathrm{is}\, \mathrm{indecomposable}\}. 
	\end{eqnarray*}
	Note that for  $\Lambda\in \widehat{T_{-1}},$  $\bar{\Lambda}:\mathbb{F}_q^\times \rightarrow\mathbb{C}^\times$ is defined  by $x\mapsto \Lambda(xI).$ 
For an indecomposable $\Lambda$, we define $\mathbb{V}_\Lambda=\mathbb{W}_{\bar{\Lambda}}\setminus \{ \chi_\Lambda^{-1}\}.$ 
	
For odd $q$, $ \widehat{\mathbb{F}_q^\times}$ is cyclic group of order $q-1$ and for each $\alpha \in \widehat{\mathbb{F}_q^\times},$ the equation $X^2=\alpha^2$ has exactly two solutions in  $\widehat{\mathbb{F}_q^\times}.$ We denote the solution distinct from  $\alpha$ by $-\alpha.$ The following theorem gives the  complete decomposition of the tensor product of any two irreducible representations of $G.$ 
	\begin{thm}\label{thm:decomp. full}
		For $\alpha,\beta,\gamma,\delta \in  \widehat{\mathbb{F}_q^\times}$
		and indecomposable $\Lambda, \Phi \in \widehat{T_{-1}}$  such that $\alpha\neq \beta$ and $\gamma\neq\delta ,$ the following hold.
		\begin{enumerate}
			\item $\chi_\alpha^1\otimes \chi_\gamma^{1} =\chi_{\alpha\gamma}^1.$
			\item $\chi_\alpha^1\otimes \chi_\gamma^{q} =\chi_{\alpha\gamma}^q.$
			\item $\chi_\alpha^1\otimes \chi_{(\gamma,\delta)}^{q+1}= \chi_{(\alpha\gamma,\alpha\delta)}^{q+1}.$
			\item $\chi_\alpha^1\otimes \chi_\Lambda^{-1} =\chi_{(\alpha\comp \det)\Lambda}^{-1}.$
			\item $\chi_\alpha^q\otimes \chi_\gamma^{q}= \chi_{\alpha\gamma}^1 \oplus \chi_{\alpha\gamma}^q  \oplus \chi_{-(\alpha\gamma)}^q  \oplus \left( \oplus_{V \in \mathbb{S}_{(\alpha\gamma)^2}} V \right) \oplus \left( \oplus_{W \in \mathbb{V}_{(\alpha\gamma)\comp \det}} W \right).$
			%\item $\chi_\alpha^q\otimes \chi_{(\gamma,\delta)}^{q+1}= \begin{cases}	\chi_{(\alpha\gamma,\alpha\delta)}^{q+1} \oplus  \chi_{\alpha\theta}^q   \oplus \chi_{-(\alpha\theta)}^q  \oplus \left( \oplus_{(\pi_1,\pi_2) \in S_{\alpha^2\gamma\delta}^\star} \chi_{(\pi_1,\pi_2)}^{q+1} \right) \oplus \left( \oplus_{\Psi \in W_{\alpha^2\gamma\delta}^\star} \chi_{\Psi}^{-1} \right)& \mathrm{if} \gamma \delta=\theta^2 \, \mathrm{for}\, \mathrm{some}\, \theta \in \widehat{\mathbb{F}_q^\times},\\  \chi_{(\alpha\gamma,\alpha\delta)}^{q+1}   \oplus \left( \oplus_{(\pi_1,\pi_2) \in S_{\alpha^2\gamma\delta}^\star} \chi_{(\pi_1,\pi_2)}^{q+1} \right) \oplus \left( \oplus_{\Psi \in W_{\alpha^2\gamma\delta}^\star} \chi_{\Psi}^{-1} \right) & \mathrm{otherwise.} \end{cases}$
			\item $\chi_\alpha^q\otimes \chi_{(\gamma,\delta)}^{q+1} = 
			\left(\oplus_{\theta \in \widehat{\mathbb{F}_q^\times}; \theta^2=\alpha^2 \gamma\delta}\chi_{\theta}^q \right)  \oplus 	2\chi_{(\alpha\gamma,\alpha\delta)}^{q+1} \oplus  \left( \oplus_{V \in \mathbb{S}_{\alpha^2\gamma\delta} \setminus\{\chi_{(\alpha\gamma,\alpha\delta)}^{q+1}\}} V \right) \oplus \left( \oplus_{W \in \mathbb{W}_{\alpha^2\gamma\delta}} W \right).$

			%\item $\chi_\alpha^q\otimes \chi_{\Lambda}^{-1}= \begin{cases} \chi_{\alpha\theta}^q \oplus \chi_{-(\alpha\theta)}^q  \oplus \left( \oplus_{(\pi_1,\pi_2) \in S_{\alpha^2\bar{\Lambda}}^\star} \chi_{(\pi_1,\pi_2)}^{q+1} \right) \oplus \left( \oplus_{\Psi \in V_{(\alpha\comp\det)\Lambda}^\star} \chi_{\Psi}^{-1} \right)&  \bar{\Lambda}=\theta^2 \, \mathrm{for}\, \mathrm{some}\, \theta \in \widehat{\mathbb{F}_q^\times},\\ \left( \oplus_{(\pi_1,\pi_2) \in S_{\alpha^2\bar{\Lambda}}^\star} \chi_{(\pi_1,\pi_2)}^{q+1} \right) \oplus \left( \oplus_{\Psi \in V_{(\alpha\comp\det)\Lambda}^\star} \chi_{\Psi}^{-1} \right) & \mathrm{otherwise.} \end{cases}$
			
			\item $\chi_\alpha^q\otimes \chi_{\Lambda}^{-1}=
			\left(\oplus_{\theta \in \widehat{\mathbb{F}_q^\times}; \theta^2=\alpha^2 \bar{\Lambda}}\chi_{\theta}^q \right)  \oplus \left( \oplus_{V \in \mathbb{S}_{\alpha^2\bar{\Lambda}}} V \right) \oplus \left( \oplus_{W \in \mathbb{V}_{(\alpha\comp\det)\Lambda}} W \right).$
			
			\item $\chi_{(\alpha,\beta)}^{q+1}\otimes\chi_{(\gamma,\delta)}^{q+1}= \left(\oplus_{\theta \in \widehat{\mathbb{F}_q^\times}; \theta^2=\alpha \beta\gamma\delta}\chi_{\theta}^q \right)  \oplus \left( \oplus_{V \in \mathbb{S}_{\alpha\beta\gamma\delta}} V \right) \oplus \left( \oplus_{W \in \mathbb{W}_{\alpha\beta\gamma\delta}} W \right) \oplus \mathrm{Ind}_B^G(\alpha\gamma ,\beta\delta)\oplus \mathrm{Ind}_B^G(\alpha\delta ,  \beta\gamma ).$

			%\item $\chi_{(\alpha, -\alpha)}^{q+1}\otimes \chi_{(\gamma, -\gamma)}^{q+1}=\chi_{\alpha\gamma}^1 \oplus \chi_{-(\alpha\gamma)}^1 \oplus 2 \chi_{\alpha\gamma}^q   \oplus 2\chi_{-(\alpha\gamma)}^q  \oplus \left( \oplus_{(\pi_1,\pi_2) \in S_{(\alpha\gamma)^2}^\star} \chi_{(\pi_1,\pi_2)}^{q+1} \right) \oplus \left( \oplus_{\Lambda \in W_{(\alpha\gamma)^2}^\star} \chi_{\Lambda}^{-1} \right)$
			\item $\chi_{(\alpha, \beta)}^{q+1}\otimes \chi_{\Lambda}^{-1} = \left(\oplus_{\theta \in \widehat{\mathbb{F}_q^\times}; \theta^2=\alpha \beta\bar{\Lambda}}\chi_{\theta}^q \right)  \oplus \left( \oplus_{V\in \mathbb{S}_{\alpha \beta\bar{\Lambda}}} V \right) \oplus \left( \oplus_{W \in \mathbb{W}_{\alpha \beta\bar{\Lambda}}} W \right) .$ 
			\item If both $\Lambda\Phi $ and $\Lambda\Phi^q$ are indecomposable, then 
			$$\chi_{\Lambda}^{-1}\otimes \chi_{\Phi}^{-1}=	\left(\oplus_{\theta \in \widehat{\mathbb{F}_q^\times}; \theta^2=\bar{\Lambda}\bar{\Phi}}\chi_{\theta}^q \right)  \oplus \left( \oplus_{V \in \mathbb{S}_{\bar{\Lambda}\bar{\Phi}}} V \right) \oplus \left( \oplus_{W \in \mathbb{W}_{\bar{\Lambda}\bar{\Phi}} \setminus \{\chi_{\Lambda\Phi}^{-1}, \chi_{\Lambda\Phi^q}^{-1}\}}W \right). $$
			\item If  $\Lambda\Phi $ is indecomposable and $\Lambda\Phi^q=\theta \comp\det$ for some $\theta \in \widehat{\mathbb{F}_q^\times},$  then 
			$$\chi_{\Lambda}^{-1}\otimes \chi_{\Phi}^{-1}=	\chi_{\theta}^1 \oplus \chi_{-\theta}^q  \oplus \left( \oplus_{V \in \mathbb{S}_{\bar{\Lambda}\bar{\Phi}}} V \right) \oplus \left( \oplus_{W \in \mathbb{W}_{\bar{\Lambda}\bar{\Phi}} \setminus \{\chi_{\Lambda\Phi}^{-1}\}} W \right) .$$
			\item If  $\Lambda\Phi=\theta \comp\det$ for some $\theta \in \widehat{\mathbb{F}_q^\times}$ and  $\Lambda\Phi^q $ is indecomposable,   then 
			$$\chi_{\Lambda}^{-1}\otimes \chi_{\Phi}^{-1}=	\chi_{\theta}^1 \oplus \chi_{-\theta}^q  \oplus \left( \oplus_{V \in \mathbb{S}_{\bar{\Lambda}\bar{\Phi}}} V \right) \oplus \left( \oplus_{W \in \mathbb{W}_{\bar{\Lambda}\bar{\Phi}} \setminus \{\chi_{\Lambda{\Phi^q}}^{-1}\}} W \right) .$$
			\item If  $\Lambda\Phi=\theta \comp\det$ and $\Lambda\Phi^q=\theta' \comp\det$  for some $\theta,\theta' \in \widehat{\mathbb{F}_q^\times},$ then $\theta'=-\theta$ and  
			$$\chi_{\Lambda}^{-1}\otimes \chi_{\Phi}^{-1}=	\chi_{\theta}^1 \oplus \chi_{-\theta}^1  \oplus \left( \oplus_{V\in \mathbb{S}_{\bar{\Lambda}\bar{\Phi}}} V \right) \oplus \left( \oplus_{W \in \mathbb{W}_{\bar{\Lambda}\bar{\Phi}} } W \right) .$$
			
		\end{enumerate}

	\end{thm}

	%Section~\ref{sec:Complete Decomposition of tensor product}. Based on the decomposition, we obtain the following classification of pairs $(V_1, V_2)$ of irreducible characters of $G$ such that $V_1 \otimes V_2$ is multiplicity free. 
	%Let $\mathrm{Irr}(G)$ denote the set of all inequivalent irreducible representations of $G.$   
	%For $r\geq 1,$ let $\mathrm{Irr}^r(G)=\{\phi \in \mathrm{Irr}(G) \mid \dim(\phi)=r\}.$
	See Section~\ref{sec:decomp. of induction from tori} for its proof.
	Let $\mathrm{Irr}(G)$ denote the set of all irreducible representations of the group $G.$ The following corollary is immediate from Theorem~\ref{thm:decomp. full} and describes the multiplicities of the irreducible constituents of tensor products:

\begin{corollary}
	Let $ V_1,V_2, W \in \mathrm{Irr}(\GL_2(\mathbb F_q)).$
	\begin{enumerate}
		\item $V_1 \otimes V_2$ is multiplicity free except for the cases 
	$\chi_\alpha^q\otimes \chi_{(\gamma,\delta)}^{q+1}$ and $\chi_{(\alpha,\beta)}^{q+1}\otimes\chi_{(\gamma,\delta)}^{q+1}.$ In both the cases, the highest multiplicity is 2 and it is due to $q$ or $(q+1)$-dimensional components. 
		\item If $W$ is a one-dimensional or a $(q-1)$-dimensional,  then $\langle  V_1 \otimes V_2, W \rangle \leq 1.$
		
	\end{enumerate}
 \end{corollary}

\begin{corollary}
	\label{thm:self dual characterisation}
	For  $V\in \mathrm{Irr}(\GL_2(\mathbb F_q)),$ $V$ is self-dual
	 %equivalent to its dual 
	 representation   if and only if $V$ belongs to one of the  following sets.
	\begin{itemize}
		\item $\{\chi_{1_{\mathbb{F}_q^\times}}^1,\chi_{-(1_{\mathbb{F}_q^\times})}^1\}.$
		\item $\{\chi_{1_{\mathbb{F}_q^\times}}^q,\chi_{-(1_{\mathbb{F}_q^\times})}^q\}.$
		\item $\{\chi_{(\alpha,\alpha^{-1})}^{q+1} \mid \alpha\in \widehat{\mathbb{F}_q^\times}, \alpha^2 \neq 1_{\mathbb{F}_q^\times}\} \cup \{\chi_{(1_{\mathbb{F}_q^\times},-(1_{\mathbb{F}_q^\times}))}^{q+1} \}.$
		\item $\{\chi_{\Lambda}^{-1} \mid \Lambda\in 
		\widehat{T_{-1}} \, \mathrm{indecomposable}, \bar{\Lambda}=1_{\mathbb{F}_q^\times} 
		\}.$
	\end{itemize}
	
\end{corollary}
It is a well known fact that if $V$ is a representation of a finite group $G,$ then  $V$ is self-dual   if and only if $V\otimes V$ contains trivial representation as a constituent.
Hence Corollary~\ref{thm:self dual characterisation} can be deduced from Theorem~\ref{thm:decomp. full}.
%\begin{corollary}
%	For $ V_1,V_2 \in \mathrm{Irr}(\GL_2(\mathbb F_q)),$ the tensor product
%		 $V_1 \otimes V_2$ is multiplicity free except for the cases $\chi_\alpha^q\otimes \chi_{(\gamma,\delta)}^{q+1}$ and $\chi_{(\alpha,\beta)}^{q+1}\otimes\chi_{(\gamma,\delta)}^{q+1}.$ In both the cases, the highest multiplicity is 2 and it is due to $q$ or $(q+1)$-dimensional components. 
%\end{corollary}

	%See Section~\ref{sec:Complete Decomposition of tensor product}  for its proof.
		For $V_1, V_2  \in \mathrm{Irr}(G),$ we say that the tensor product $V_1\otimes V_2$ 
		has {\it unique decomposition property}, if whenever there exist $ W_1, W_2 \in \mathrm{Irr}(G)$  such that $V_1\otimes V_2= W_1\otimes W_2,$ then
		$V_i=\chi_i\otimes W_{\sigma(i)} ;\,i=1,2$, for some  one-dimensional representations $\chi_1$ and $\chi_2$ of $G$  and some permutation $\sigma$ of the set $\{1,2\}.$ This problem has been considered in \cite{MR2123935} for simple algebraic groups and in \cite{MR2980495} for Kac-Moody algebras.
		%$\{W_1,W_2\}=\{\chi\otimes V_1, \phi \otimes V_2\}$ for some  one-dimensional representations $\chi$ and $\phi$ of $G.$ 
	We consider the unique decomposition problem of $\GL_2{(\mathbb{F}_q)}$ and obtain the following result: 
	\begin{thm} \label{thm:unique decomposition}
		For  $V_1,V_2 \in \mathrm{Irr}(\GL_2(\mathbb F_q))$ such that $\{\dim(V_1),\dim(V_2)\}\neq \{q-1,q+1\},$ the tensor product  $V_1\otimes V_2$  has unique decomposition property.
	\end{thm}
 The following theorem 
 %(see Section~\ref{sec:Unique Decomposition of tensor product} for its proof) 
 addresses the remaining case:
	\begin{thm}\label{thm:unique decomp q-1 and q+1}
		Let $\Lambda, \Phi \in \widehat{T_{-1}}$  be two indecomposable characters and  $\alpha, \beta, \gamma, \delta\in \widehat{\mathbb{F}_q^\times} $ be such that $\alpha\neq\beta$ and $ \gamma\neq  \delta.$  Then   $\chi_{(\alpha,\beta)}^{q+1}\otimes\chi_{\Lambda}^{-1}=\chi_{(\gamma,\delta)}^{q+1}\otimes\chi_{\Phi}^{-1}$ if and only if $\alpha \beta \bar{\Lambda}=\gamma \delta \bar{\Phi}.$ In particular, for $q\neq 3,$ the tensor product $\chi_{(\alpha,\beta)}^{q+1}\otimes\chi_{\Lambda}^{-1}$
		does not have unique decomposition property.
		
	\end{thm}

	\begin{remark}
		For $q=3,$ there are exactly  three $(q-1)$-dimensional  representations and  one $(q+1)$-dimensional representation in $\mathrm{Irr}(G).$ It is easy show that the tensor product of those representations has unique decomposition property.
	\end{remark}
	See Section~\ref{sec:Unique Decomposition of tensor product} for proofs of Theorems~\ref{thm:unique decomposition} and \ref{thm:unique decomp q-1 and q+1}.

 \subsection*{Acknowledgement:}
 The authors would like to express their deepest gratitude to Pooja Singla for careful scrutiny of the draft and many helpful conversations.
 \section{Proof of Theorem \ref{thm:decomp. full}}\label{sec:decomp. of induction from tori}
 We first  introduce some notations. For an irreducible representation $\phi$ of a subgroup $H$ of $G,$ we denote the set of all inequivalent irreducible constituents of $\mathrm{Ind}_H^G\phi$ by $\mathrm{Irr}(G\mid\phi).$
 	Let $B, T_1,T_{-1}, Z$ and $U$ be the standard borel subgroup consisting of upper-triangular matrices, isotropic torus, anisotropic torus, center and standard unipotent subgroup of $G$ respectively. 
 	%For  $\Lambda \in \widehat{T_{-1}},$ we define a character $\bar{\Lambda}:\mathbb{F}_q^\times \rightarrow\mathbb{C}^\times$ by $x\mapsto \Lambda(xI).$ 
 A character $\Lambda \in \widehat{T_{-1}}$ is called {\it decomposable} if there exists $\alpha \in \widehat{\mathbb{F}_q^\times}$ such that   $\Lambda=\alpha \comp \det .$  	A character $\Lambda \in \widehat{T_{-1}}$ is called {\it indecomposable} if it is not decomposable.
 For  $ \alpha, \beta \in \widehat{\mathbb{F}_q^\times},$ define a one-dimensional representation $(\alpha, \beta): B \rightarrow \mathbb{C}^\times,$ by $\mat{x}{z}{0}{y}\mapsto \alpha(x) \beta(y) .$ We  
use the same notation  $(\alpha, \beta)$ for the restriction of $(\alpha, \beta)$ to $T_1.$ Fix a non-trivial character $\psi$ of $U.$ For a character $\rho\in \widehat{T_{i} }, i\in\{\pm1\}$, define a character  $\rho\psi: ZU \rightarrow \mathbb{C}^\times$ by $\mat{x}{y}{0}{x}\mapsto \rho(\mat{x}{0}{0}{x})\psi(\mat{1}{x^{-1}y}{0}{1}) .$ 

 To prove Theorem~\ref{thm:decomp. full}, we use the following result:

\begin{thm} \cite[Theorem~3.1]{MR1757476}
	\label{thm:Jose Pantoja tensor product}
	For $\alpha,\beta,\gamma,\delta \in  \widehat{\mathbb{F}_q^\times}$
	and indecomposable $\Lambda, \Phi \in \widehat{T_{-1}}$ such that $\alpha\neq \beta$ and $\gamma\neq\delta ,$ the following hold.
	\begin{enumerate}
		\item$\chi_\alpha^q\otimes \chi_\gamma^q=\Ind_{T_{-1}}^G((\alpha\gamma)\comp \det)+\chi_{\alpha\gamma}^q.$ 
		\item $\chi_\alpha^q\otimes \chi_{( \gamma,\delta)}^{q+1}=\Ind_{T_1}^G(\alpha\gamma, \alpha\delta).$ 
		\item $\chi_\alpha^q\otimes \chi_{\Lambda}^{-1}=\Ind_{T_{-1}}^G(\alpha\comp \det)\Lambda.$
		\item $\chi_{(\alpha, \beta)}^{q+1}\otimes \chi_{(\gamma, \delta)}^{q+1}=\begin{cases}
			\Ind_{T_1}^G(\alpha\gamma, \beta\delta)+\chi_{(\beta\gamma,\alpha\delta)}^{q+1} & \mathrm{if}\, \beta\gamma\neq\alpha\delta,\\
			\Ind_{T_1}^G(\alpha\gamma, \beta\delta)+\chi_{\beta\gamma}^1+\chi_{\beta\gamma}^q & \mathrm{if}\, \beta\gamma=\alpha\delta.
		\end{cases}$
		\item $\chi_{(\alpha, \beta)}^{q+1}\otimes \chi_{\Lambda}^{-1}=
		\begin{cases}
			\Ind_{T_1}^G(\alpha\beta,\bar{\Lambda})-\chi_{(\alpha\beta,\bar{\Lambda})}^{q+1} & \mathrm{if}\, \alpha\beta\neq \bar{\Lambda},\\
			\Ind_{T_1}^G(\alpha\beta,\bar{\Lambda})-\chi_{\alpha\beta}^1 -\chi_{\alpha\beta}^q & \mathrm{if}\, \alpha\beta= \bar{\Lambda}.
		\end{cases}$ 
		\item $\chi_\Lambda^{-1}\otimes\chi_\Phi^{-1}=\begin{cases}
			\Ind_{T_{-1}}^G\Lambda\Phi^q-\chi_{\Lambda\Phi}^{-1} & \mathrm{if} \, \Lambda\Phi  \, \mathrm{is}\, \mathrm{indecomposable,} \\
			\Ind_{T_{-1}}^G\Lambda\Phi^q-\chi_{\theta}^{q}+ \chi_{\theta}^{1}  & \mathrm{if} \, \Lambda\Phi = \theta \comp \det \, \mathrm{for}\, \mathrm{some}\, \theta \in \widehat{\mathbb{F}_q^\times}.
		\end{cases} $ 
		
	\end{enumerate}
\end{thm}
Thus, to prove Theorem~\ref{thm:decomp. full}, it is enough to decompose the representations $\mathrm{Ind}_{T_{-1}}^G\Lambda$ and $\mathrm{Ind}_{T_{1}}^G(\alpha,\beta).$  Next lemma shows that decomposing $\Ind_{ZU}^G\rho\psi$ would suffice.

\begin{lemma}
	\label{lem:torus in terms of ZU to G}
	For $\Lambda\in\widehat{T_{-1}}$ and $(\alpha,\beta)\in\widehat{T_1},$ the following hold.
	\begin{enumerate}
		\item $\mathrm{Ind}_{T_{-1}}^G\Lambda=\begin{cases}
			\mathrm{Ind}_{ZU}^G\Lambda\psi-\chi_{\Lambda}^{-1}&   \mathrm{if}\,\Lambda  \, \mathrm{is}\, \mathrm{indecomposable,} \\  \mathrm{Ind}_{ZU}^G\Lambda\psi-\chi_\alpha^q+\chi_\alpha^1 &\mathrm{if}\, \Lambda=\alpha\comp\det\, \mathrm{for}\, \mathrm{some}\, \alpha \in \widehat{\mathbb{F}_q^\times}.
		\end{cases}$
		\item  $\mathrm{Ind}_{T_{1}}^G(\alpha,\beta)=\begin{cases}
			\mathrm{Ind}_{ZU}^G(\alpha,\beta)\psi+\chi_{(\alpha,\beta)}^{q+1}&   \mathrm{if}\,\alpha\neq\beta, \\  \mathrm{Ind}_{ZU}^G(\alpha,\alpha)\psi+\chi_\alpha^q+\chi_\alpha^1 &\mathrm{if}\, \alpha =\beta.
		\end{cases}$
		
	\end{enumerate}
\end{lemma}
\begin{proof}
	The lemma directly follows from \cite[Proposition~2.1 and 2.2]{MR1757476}.
\end{proof}

%By the above lemma, for proving theorems \ref{thm:T-1 to G decomposition} and \ref{thm:T1 to G decomposition}, it is enough to have decomposition of $\mathrm{Ind}_{ZU}^G\rho\psi$, where $\rho\in\widehat{Z}$ and $\psi\in\widehat{U}$, $\psi$ non trivial. This is given by the following results:

\begin{lemma}\label{lem:ZU to G}
	For $\rho\in \widehat{Z}$ and a non-trivial character  $ \psi$ of $U,$ 
	$$\mathrm{Ind}_{ZU}^G(\rho\psi)=	\left(\oplus_{\alpha \in \widehat{\mathbb{F}_q^\times}; \alpha^2=\rho}\chi_{\alpha}^q \right)  \oplus \left( \oplus_{V \in \mathbb{S}_{\rho}} V \right) \oplus \left( \oplus_{W \in \mathbb{W}_{\rho}} W \right) .$$
\end{lemma}

\begin{proof}
	By \cite[Theorem 16.1]{MR2270898}, we have that $  \mathrm{Ind}_{U}^G\psi$ is multiplicity free and  $  \mathrm{Irr}(G \mid \psi)= \{\pi\in \mathrm{Irr}(G) \mid  \dim(\pi)\neq 1\}.$
	Also, it is easy to see that $ \mathrm{Irr}(ZU \mid \psi)= \{\rho\psi \mid  \rho \in \mathrm{Irr}(Z)\}.$ Therefore, for each $ \rho \in \mathrm{Irr}(Z), $ the representation $\mathrm{Ind}_{ZU}^G(\rho\psi)$ is multiplicity free and  $  \mathrm{Irr}(G \mid \rho\psi)= \{\pi\in \mathrm{Irr}(G) \mid  \dim(\pi)\neq 1, \langle \mathrm{Res}^{G}_Z\pi, \rho \rangle \neq 0 \}.$ Hence the lemma follows from the well known facts that 
	$\mathrm{Res}^{G}_Z\chi_\alpha^q=q \alpha^2,$ $\mathrm{Res}^{G}_Z \chi_{(\alpha,\beta)}^{q+1}=(q+1)\alpha\beta$ and $\mathrm{Res}^{G}_Z \chi_{\Lambda}^{-1}=(q-1)\bar{\Lambda}.$
	
\end{proof}	
Plugging the decomposition of $\Ind_{ZU}^G(\rho\psi)$ in lemma \ref{lem:torus in terms of ZU to G} gives the decomposition of induction from tori:

\begin{thm}
	\label{thm:T-1 to G decomposition}
	For  $\Lambda \in \widehat{T_{-1}},$
	%, and suppose $\Lambda^2(y)=\zeta_{q-1}^t$, for some $1\leq t\leq (q-1)$ where $y$ is a generator of $K^\times$ then 
	the decomposition of $\mathrm{Ind}_{T_{-1}}^G\Lambda$ is as follows. 
	\begin{enumerate} 
		\item If $\Lambda$ is a decomposable character with  $\Lambda=\beta \comp \det $ for some $\beta \in \widehat{\mathbb{F}_q^\times},$  then 
		\begin{equation*}
			\mathrm{Ind}_{T_{-1}}^G\Lambda =  \chi_{\beta}^1 \oplus \chi_{-\beta}^q  \oplus \left( \oplus_{V \in \mathbb{S}_{\bar{\Lambda}}} V \right) \oplus \left( \oplus_{W \in \mathbb{W}_\Lambda} W \right).
		\end{equation*}
		
		\item If $\Lambda$ is a indecomposable character such that $\bar{\Lambda}=\beta^2$ for some  $\beta\in \widehat{\mathbb{F}_q^\times},$ then
		\begin{equation*}
			\mathrm{Ind}_{T_{-1}}^G\Lambda =\chi_{\beta}^q \oplus \chi_{-\beta}^q  \oplus \left( \oplus_{V \in \mathbb{S}_{\bar{\Lambda}}} V \right) \oplus \left( \oplus_{W \in \mathbb{V}_\Lambda} W \right). 
		\end{equation*}

		\item If $\Lambda$ is a indecomposable character such that $\bar{\Lambda}\neq\beta^2$ for any $\beta\in \widehat{\mathbb{F}_q^\times},$ then
		$$ \mathrm{Ind}_{T_{-1}}^G\Lambda =\left( \oplus_{V \in \mathbb{S}_{\bar{\Lambda}}} V \right) \oplus \left( \oplus_{W \in \mathbb{V}_\Lambda} W \right). $$ 
	\end{enumerate}
	
\end{thm}

\begin{thm}
	\label{thm:T1 to G decomposition}
	For  $\alpha,\beta \in \widehat{\mathbb{F}_q^\times},$
	%, and suppose $\alpha^2(y)=\zeta_{q-1}^t$, for some $1\leq t\leq (q-1)$ where $y$ is a generator of $K^\times$ then 
	the decomposition of $\mathrm{Ind}_{T_{1}}^G(\alpha, \beta)$ is as follows. 
	\begin{enumerate} 
		\item If $\alpha= \beta,$   then 
		\begin{equation*}
			\mathrm{Ind}_{T_{1}}^G(\alpha, \alpha) = \chi_\alpha^1 \oplus 2 \chi_\alpha^q   \oplus \chi_{-\alpha}^q  \oplus \left( \oplus_{V \in \mathbb{S}_{\alpha^2}} V \right) \oplus \left( \oplus_{W \in \mathbb{W}_{\alpha^2}} W \right).
		\end{equation*}
		
		\item If $\alpha\neq \beta$ and  $\alpha\beta=\gamma^2$ for some  $\gamma\in \widehat{\mathbb{F}_q^\times},$ then
		\begin{equation*}
			\mathrm{Ind}_{T_{1}}^G(\alpha, \beta)=  \chi_{(\alpha,\beta)}^{q+1} \oplus  \chi_\gamma^q   \oplus \chi_{-\gamma}^q  \oplus \left( \oplus_{V\in \mathbb{S}_{\alpha\beta}}V \right) \oplus \left( \oplus_{W \in \mathbb{W}_{\alpha\beta}}W \right).
		\end{equation*}
		
		\item If $\alpha\neq \beta$ and  $\alpha\beta\neq\gamma^2$ for all  $\gamma\in \widehat{\mathbb{F}_q^\times},$ then
		\begin{equation*}
			\mathrm{Ind}_{T_{1}}^G(\alpha, \beta) =  \chi_{(\alpha,\beta)}^{q+1} \oplus  \left( \oplus_{V\in \mathbb{S}_{\alpha\beta}} V \right) \oplus \left( \oplus_{W \in \mathbb{W}_{\alpha\beta}} W \right).
		\end{equation*}

	\end{enumerate}
	
\end{thm}
\begin{remark*}

	         The following observations are evident from the above theorems,
	\begin{enumerate}
			\item  $\mathrm{Ind}_{T_{-1}}^G\Lambda$ is multiplicity free.
		\item For $\alpha\neq\beta$, since $\chi_{(\alpha,\beta)}^{q+1}\in \mathbb{S}_{\alpha\beta}$, the multiplicity of $ \chi_{(\alpha,\beta)}^{q+1}$  in the decomposition of $\mathrm{Ind}_{T_{1}}^G(\alpha, \beta)$ is two.
	
	\end{enumerate}

\end{remark*}

\begin{lemma}\label{lem:theta'=-theta}
	For indecomposable $\Lambda, \Phi \in \widehat{T_{-1}},$ if  $\Lambda\Phi=\theta \comp\det$ and $\Lambda\Phi^q=\theta' \comp\det$  for some $\theta,\theta' \in \widehat{\mathbb{F}_q^\times},$ then $\theta'=-\theta.$ 
\end{lemma}
\begin{proof}
	Note that $\theta^2=\overline{(\Lambda\Phi)}=\bar{\Lambda}\bar{\Phi}=\bar{\Lambda}\bar{\Phi}^q=\overline{(\Lambda\Phi^q)}=\theta^{\prime 2}.$ Therefore $\theta'\in \{\theta, -\theta\}.$ If $\theta' =\theta,$ then $\Lambda\Phi^q=\Lambda\Phi$ and hence $\Phi^{q-1}=1_{T_{-1}}.$ Therefore, since  $\widehat{T_{-1}}$ is a cyclic group of order $q^2-1$  (because $\widehat{T_{-1}}\cong T_{-1}\cong \mathbb{F}_{q^2}^\times$), $\Phi$ must be in the subgroup of $\widehat{T_{-1}}$ with order $q-1.$ But $\{\alpha\comp\det \mid \alpha \in \widehat{\mathbb{F}_q^\times}\}$ is the subgroup of $\widehat{T_{-1}}$ with order $q-1.$ Therefore $\Phi$ is decomposable, which is a contradiction. Hence we must have $\theta'=-\theta.$ 
\end{proof}
\begin{proof}[\bf Proof of Theorem~\ref{thm:decomp. full}]
	%\footnote{add proof}
	Note that (1)-(4) of Theorem~\ref{thm:decomp. full}, directly follows from the character table of $G$ (see \cite{MR2319497}). Similarly, (5)-(13) of Theorem~\ref{thm:decomp. full} directly follow from Theorems~\ref{thm:Jose Pantoja tensor product}, \ref{thm:T-1 to G decomposition} and \ref{thm:T1 to G decomposition},
	Lemma~\ref{lem:theta'=-theta} and the fact that, for $\theta\in \widehat{\mathbb{F}_q^\times}$ and indecomposable $\Psi \in \widehat{T_{-1}},$ the character $(\theta\comp \det)\Psi$ is indecomposable.
\end{proof}

\section{Unique decomposition property of the tensor product % of irreducible representations of $\GL_2(k)$
}\label{sec:Unique Decomposition of tensor product}

In this section, we prove Theorems~\ref{thm:unique decomposition} and  \ref{thm:unique decomp q-1 and q+1}. For that we need the following lemmas.

\begin{lemma}\label{lem:unique decomp.}
	For indecomposable $\Lambda,\Lambda'\in \widehat{T_{-1}},$ $\mathbb{V}_\Lambda=\mathbb{V}_{\Lambda'}$ if and only if $\chi_\Lambda^{-1}= \chi_{\Lambda'}^{-1}.$

\end{lemma}
\begin{proof} 
	If $\chi_\Lambda^{-1}= \chi_{\Lambda'}^{-1},$ then  $\bar{\Lambda}=\bar{\Lambda'}.$ Therefore by definition of $\mathbb{V}_\Lambda,$ 
	%and $\mathbb{W}_{\bar{\Lambda}},$ 
	it is clear that $\mathbb{V}_\Lambda=\mathbb{V}_\Lambda'.$ To show the converse, suppose  $\mathbb{V}_\Lambda=\mathbb{V}_{\Lambda'}.$  
	We first show  that   $\bar{\Lambda'}=\bar{\Lambda}.$ If $\mathbb{V}_\Lambda \neq \emptyset,$ choose $\chi_\Phi^{-1}\in \mathbb{V}_\Lambda$ and by definition of $\mathbb{V}_\Lambda$ and $\mathbb{V}_{\Lambda'},$  we must have  $\bar{\Lambda}=\bar{\Phi}=\bar{\Lambda'}.$ If $\mathbb{V}_\Lambda = \emptyset,$ then $|\mathbb{V}_\Lambda |=0.$ From the definition of $\mathbb{V}_\Lambda ,$ it is easy to see that 
	$$ |\mathbb{V}_\Lambda|=\begin{cases}
		\frac{q-3}{2} &\mathrm{if}\, \bar{\Lambda}=\gamma^2 \, \mathrm{for}\, \mathrm{some}\, \gamma \in \widehat{\mathbb{F}_q^\times},\\
		\frac{q-1}{2} &\mathrm{otherwise.} \end{cases}$$
	Therefore  we must have   $q=3$ and  $\{\bar{\Lambda}, \bar{\Lambda'} \} \subseteq \{\gamma^2 \mid \gamma \in \widehat{\mathbb{F}_q^\times}\}.$ Since $\{\gamma^2 \mid \gamma \in \widehat{\mathbb{F}_3^\times}\}=\{1_{\mathbb{F}_3^\times}\},$ we obtain $\bar{\Lambda}=1_{\mathbb{F}_3^\times} =\bar{\Lambda'}.$
	Next we show that $\chi_\Lambda^{-1}= \chi_{\Lambda'}^{-1}.$
	Note that  $\bar{\Lambda} =\bar{\Lambda'}$ implies $\mathbb{W}_{\bar{\Lambda}}=\mathbb{W}_{\bar{\Lambda'}}.$ Therefore
	$$\{\chi_\Lambda^{-1}\}=  \mathbb{W}_{\bar{\Lambda}} \setminus   \mathbb{V}_{\Lambda}= \mathbb{W}_{\bar{\Lambda'}} \setminus \mathbb{V}_{\Lambda'} = \{\chi_{\Lambda'}^{-1}\}.$$
	Hence the lemma holds.
\end{proof}

%\begin{lemma}\label{lem:uniqe decomp2}
%	For $\alpha,\beta,\gamma,\delta, \alpha',\beta',\gamma',\delta' \in \widehat{\mathbb{F}_q^\times},$
%	$ \mathrm{Ind}_B^G(\alpha,\beta)\oplus \mathrm{Ind}_B^G(\gamma , \delta)=\mathrm{Ind}_B^G(\alpha',\beta')\oplus \mathrm{Ind}_B^G(\gamma' , \delta')$ if and only if $ \{\{\alpha,\beta\}, \{\gamma , \delta\}\}=\{\{\alpha',\beta'\}, \{\gamma' , \delta'\}\}$
%\end{lemma}
%\begin{proof}
%	For $\alpha,\beta \in \widehat{\mathbb{F}_q^\times},$ we have $$\mathrm{Ind}_B^G(\alpha,\beta)=\begin{cases}
%		\chi_\alpha^1 \oplus\chi_\alpha^q & \mathrm{if}\, \alpha=\beta,\\
%		\chi_{\alpha,\beta}^{q+1}& \mathrm{if}\, \alpha\neq\beta.
%	\end{cases}$$
%	Therefore it is easy to see that 	$ \mathrm{Ind}_B^G(\alpha,\beta)\oplus \mathrm{Ind}_B^G(\gamma , \delta)=\mathrm{Ind}_B^G(\alpha',\beta')\oplus \mathrm{Ind}_B^G(\gamma' , \delta')$ if and only if 
%	$ \{\mathrm{Ind}_B^G(\alpha,\beta),  \mathrm{Ind}_B^G(\gamma , \delta)\}= \{\mathrm{Ind}_B^G(\alpha',\beta'), \mathrm{Ind}_B^G(\gamma' , \delta') \},$ which is if and only if 
%	$ \{\{\alpha,\beta\}, \{\gamma , \delta\}\}=\{\{\alpha',\beta'\}, \{\gamma' , \delta'\}\}.$
%\end{proof}

\begin{lemma}\label{lem:cuspidal with cuspidal}
	Let $\Lambda,\Phi,\Lambda', \Phi'\in \widehat{T_{-1}}$ be indecomposable characters. If 
	$\chi_{\Lambda}^{-1}\otimes \chi_{\Phi}^{-1}=\chi_{\Lambda'}^{-1}\otimes \chi_{\Phi'}^{-1},$ then 
	either
	$ \Lambda'\Phi'\in \{\Lambda\Phi, \Lambda^q\Phi^q\}\, \mathrm{and}\,\Lambda'\Phi^{\prime q} \in \{\Lambda\Phi^q, \Lambda^q\Phi\},  \, \mathrm{or} \, \Lambda'\Phi^{\prime q} \in \{\Lambda\Phi, \Lambda^q\Phi^q\}\, \mathrm{and}\, \Lambda'\Phi'  \in \{\Lambda\Phi^q, \Lambda^q\Phi\}. $

	% $ \{\{\Lambda\Phi, \Lambda^q\Phi^q\},\{\Lambda\Phi^q, \Lambda^q\Phi\}\}=\{\{\Lambda'\Phi', \Lambda^{\prime q}\Phi^{\prime q}\},\{\Lambda'\Phi^{\prime q}, \Lambda^{\prime q}\Phi'\}\}.$
\end{lemma}
\begin{proof}
	For $\Psi\in \widehat{T_{-1}},$ let 
	$$f(\Psi)=\begin{cases}
		\chi_\alpha^q -\chi_\alpha^1 & \mathrm{if} \, \Psi=\alpha\comp \det  \, \mathrm{for}\, \mathrm{some}\, \alpha \in \widehat{\mathbb{F}_q^\times},\\
		\chi_\Psi^{-1} & \mathrm{if} \, \Psi \, \mathrm{is} \, \mathrm{indecomposable.}
	\end{cases}$$
	Note that for indecomposable $\Lambda,\Phi\in \widehat{T_{-1}},$ by Theorem~\ref{thm:decomp. full}(10)-(13), we have
	\begin{equation}\label{eqn:asw}
		(\chi_{\Lambda}^{-1}\otimes \chi_{\Phi}^{-1})\oplus f(\Lambda\Phi)\oplus f(\Lambda\Phi^q)= 	\left(\oplus_{\alpha \in \widehat{\mathbb{F}_q^\times}; \alpha^2=\bar{\Lambda}\bar{\Phi}}\chi_{\alpha}^q \right)  \oplus \left( \oplus_{\pi \in \mathbb{S}_{\bar{\Lambda}\bar{\Phi}}} \pi \right) \oplus \left( \oplus_{\rho \in \mathbb{W}_{\bar{\Lambda}\bar{\Phi}}} \rho\right) .
	\end{equation}

	Let $\Lambda,\Phi,\Lambda', \Phi'\in \widehat{T_{-1}}$ be indecomposable characters such that
	$\chi_{\Lambda}^{-1}\otimes \chi_{\Phi}^{-1}=\chi_{\Lambda'}^{-1}\otimes \chi_{\Phi'}^{-1}.$ 
	By restricting both sides to $Z,$ we obtain $\bar{\Lambda}\bar{\Phi}=\bar{\Lambda'}\bar{\Phi'}.$
	Therefore by (\ref{eqn:asw}), we must have  
	$$(\chi_{\Lambda}^{-1}\otimes \chi_{\Phi}^{-1})\oplus f(\Lambda\Phi)\oplus f(\Lambda\Phi^q) = (\chi_{\Lambda'}^{-1}\otimes \chi_{\Phi'}^{-1})\oplus f(\Lambda'\Phi')\oplus f(\Lambda'\Phi^{\prime q}),$$
	which implies $f(\Lambda\Phi)\oplus f(\Lambda\Phi^q) = f(\Lambda'\Phi')\oplus f(\Lambda'\Phi^{\prime q}).$
	Hence, by definition of $f(\cdot),$ either $ f(\Lambda\Phi)= f(\Lambda'\Phi')$ and   $ f(\Lambda\Phi^q) =  f(\Lambda'\Phi^{\prime q}),$ or 
	$f(\Lambda\Phi)=  f(\Lambda'\Phi^{\prime q})$ and $ f(\Lambda\Phi^q) = f(\Lambda'\Phi') . $ 
	Note that for $\Psi_1, \Psi_2 \in \widehat{T_{-1}},$ $f(\Psi_1)= f(\Psi_2)$ if and only if $\Psi_2 \in \{\Psi_1,\Psi_1^q\}.$ 
	Therefore we obtain that either	$ \Lambda'\Phi'\in \{\Lambda\Phi, \Lambda^q\Phi^q\}\, \mathrm{and}\,\\ \Lambda'\Phi^{\prime q}
	 \in \{\Lambda\Phi^q, \Lambda^q\Phi\},  \, \mathrm{or} \, \Lambda'\Phi^{\prime q} \in \{\Lambda\Phi, \Lambda^q\Phi^q\}\, \mathrm{and}\, \Lambda'\Phi'  \in \{\Lambda\Phi^q, \Lambda^q\Phi\}. $
\end{proof}

%\begin{thm}
%	For  $V_1,V_2 \in \mathrm{Irr}(G)$ a such that $\{\dim(V_1),\dim(V_2)\}\neq \{q-1,q+1\},$ the tensor product  $V_1\otimes V_2$ has unique decomposition property.
%\end{thm}
\begin{proof}[\bf Proof of Theorem~\ref{thm:unique decomposition}] 
	Let $V_1,V_2,W_1,W_2 \in \mathrm{Irr}(G)$ be such that  $V_1\otimes V_2= W_1\otimes W_2.$ Since $\dim(V)\in \{1,q-1,q,q+1\}$ for all $V  \in \mathrm{Irr}(G),$ we must have $\{\dim(V_1),\dim(V_2)\}=\{\dim(W_1),\dim(W_2)\}.$ We assume $\dim(V_1)>1$ and  $\dim(V_2)>1$ because otherwise the unique decomposition follows trivially. 
	We first observe that for any pair $\chi_\alpha^q ,\chi_\beta^q \in \mathrm{Irr}^q(G),$ 
	%by Theorem~\ref{thm:decomp. full}(2),  
	we have $\chi_\beta^q = \chi_{\alpha^{-1}\beta}^1 \otimes \chi_\alpha^q. $ Therefore if $\dim(V_1)=q$ and $\dim(V_2)=q,$ then unique decomposition follows trivially.

	\noindent
	{\bf \underline{Case 1}} $(\dim(V_1)=q$ and $\dim(V_2)=q+1):$
	Suppose $\chi_{\alpha}^q\otimes\chi_{(\beta,\gamma)}^{q+1}=    \chi_{\alpha'}^q\otimes\chi_{(\beta',\gamma')}^{q+1}.$
	By equating the $(q+1)$-dimensional irreducible constituents with multiplicity two (from Theorem~\ref{thm:decomp. full}(6)) in both sides, we obtain $\chi_{(\alpha\beta,\alpha\gamma)}^{q+1}=    \chi_{(\alpha'\beta',\alpha'\gamma')}^{q+1}.$ Therefore 
	% by Theorem~\ref{thm:decomp. full}(3), we obtain
	$  \chi_{(\beta',\gamma')}^{q+1}=  \chi_{\alpha^{\prime-1}\alpha}^1 \otimes \chi_{(\beta,\gamma)}^{q+1} .$

	\noindent
	{\bf \underline{Case 2}} $(\dim(V_1)=q$ and $\dim(V_2)=q-1):$
	Suppose $\chi_{\alpha}^q\otimes\chi_{\Lambda}^{-1}=    \chi_{\alpha'}^q\otimes\chi_{\Lambda'}^{-1}.$
	By equating the set of  $(q-1)$-dimensional irreducible constituents   in both sides (from Theorem~\ref{thm:decomp. full}(7)), we obtain that $\mathbb{V}_{(\alpha\comp\det) \Lambda}=\mathbb{V}_{(\alpha'\comp\det)\Lambda'}.$ Since both $(\alpha\comp\det) \Lambda$ and $(\alpha'\comp\det) \Lambda'$ are indecomposable, by Lemma~\ref{lem:unique decomp.}, we must have 
	$\chi_{(\alpha\comp\det)\Lambda}^{-1}= \chi_{(\alpha'\comp\det)\Lambda'}^{-1}.$ %This together with Theorem~\ref{thm:decomp. full}(4) imply that 
	Therefore
	$ \chi_{\Lambda'}^{-1}= \chi_{\alpha^{\prime-1}\alpha}^1 \otimes \chi_{\Lambda}^{-1}.$
	%	Also, by Lemma~\ref{lem:unique decomp.}(1), we have $\chi_{\alpha'}^q=\chi_{\alpha^{-1}\alpha^{\prime}}^1 \otimes\chi_{\alpha}^q$

	\noindent
	{\bf \underline{Case 3}} $(\dim(V_1)=q+1$ and $\dim(V_2)=q+1):$
	Suppose $\chi_{(\alpha,\beta)}^{q+1}\otimes\chi_{(\gamma, \delta)}^{q+1}=    \chi_{(\alpha',\beta')}^{q+1}\otimes\chi_{(\gamma', \delta')}^{q+1}.$
	%By equating the set of  $(q-1)$-dimensional irreducible constituents   in both sides (from Theorem~\ref{thm:decomp. full}(8)), we obtain that $W_{\alpha\beta\gamma\delta}^\star=W_{\alpha'\beta'\gamma'\delta'}^\star.$ Therefore by Lemma~\ref{lem:unique decomp.}(1), we have $\alpha\beta\gamma\delta=\alpha'\beta'\gamma'\delta'.$ 
	By restricting both sides to $Z,$ we obtain $\alpha\beta\gamma\delta=\alpha'\beta'\gamma'\delta'.$  
	%Hence $S_{\alpha\beta\gamma\delta}^\star=S_{\alpha'\beta'\gamma'\delta'}^\star.$ 
	Therefore by Theorem~\ref{thm:decomp. full}(8), we must have 
	$$ \mathrm{Ind}_B^G(\alpha\gamma ,\beta\delta)\oplus \mathrm{Ind}_B^G(\alpha\delta ,  \beta\gamma )=\mathrm{Ind}_B^G(\alpha'\gamma' ,\beta'\delta')\oplus \mathrm{Ind}_B^G(\alpha'\delta' ,  \beta'\gamma' ).$$
	This implies $ \{\mathrm{Ind}_B^G(\alpha\gamma ,\beta\delta), \mathrm{Ind}_B^G(\alpha\delta ,  \beta\gamma )\}= \{ \mathrm{Ind}_B^G(\alpha'\gamma' ,\beta'\delta'), \mathrm{Ind}_B^G(\alpha'\delta' ,  \beta'\gamma' )\},$ which is equivalent to 
	$\{\{\alpha\gamma ,\beta\delta\},\{\alpha\delta ,  \beta\gamma\}\}=\{\{\alpha'\gamma' ,\beta'\delta'\},\{\alpha'\delta' ,  \beta'\gamma'\}\}.$ So there are eight  possible cases. For example,  $\alpha\gamma=\alpha'\gamma',$ $\beta\delta=\beta'\delta',$ $\alpha\delta=\alpha'\delta' $ and $ \beta\gamma=\beta'\gamma'.$
	For  this case,  we have 
	$(\gamma^{\prime-1}\gamma)\alpha=\alpha'=(\delta^{\prime-1}\delta)\alpha$ and  $(\gamma^{\prime-1}\gamma)\beta=\beta'=(\delta^{\prime-1}\delta)\beta.$ Therefore we must have $\gamma^{\prime-1}\gamma=\delta^{\prime-1}\delta$ and we obtain 
	$(\alpha', \beta')=(\gamma^{\prime-1}\gamma\alpha, \gamma^{\prime-1}\gamma\beta)$  and $(\gamma', \delta')=(\gamma^{-1}\gamma'\gamma, \gamma^{-1}\gamma'\delta).$
	Hence 
	%by  Theorem~\ref{thm:decomp. full}(3), we get that 
	$$\chi_{(\alpha',\beta')}^{q+1}= \chi_{\gamma^{\prime-1}\gamma}^1\otimes\chi_{(\alpha,\beta)}^{q+1} \, \mathrm{and}\, \chi_{(\gamma',\delta')}^{q+1} = \chi_{\gamma^{-1}\gamma'}^1   \otimes\chi_{(\gamma, \delta)}^{q+1}. $$
	By similar arguments,  we can show the unique decomposition property for other cases.

	\noindent
	{\bf \underline{Case 4}} $(\dim(V_1)=q-1$ and $\dim(V_2)=q-1):$
	Suppose 
	$\chi_{\Lambda}^{-1}\otimes \chi_{\Phi}^{-1}=\chi_{\Lambda'}^{-1}\otimes \chi_{\Phi'}^{-1}.$
	By Lemma~\ref{lem:cuspidal with cuspidal}, we have 
	either
	$ \Lambda'\Phi'\in \{\Lambda\Phi, \Lambda^q\Phi^q\}\, \mathrm{and}\,\Lambda'\Phi^{\prime q} \in \{\Lambda\Phi^q, \Lambda^q\Phi\},  \, \mathrm{or} \, \Lambda'\Phi^{\prime q} \in \{\Lambda\Phi, \Lambda^q\Phi^q\}\, \mathrm{and}\, \Lambda'\Phi'  \in \{\Lambda\Phi^q, \Lambda^q\Phi\}. $
	So there are eight  possible cases.  For example,  consider the following two cases.
	\begin{enumerate}
		\item $ \Lambda'\Phi^{\prime q}= \Lambda\Phi^{q}  \, \mathrm{and}\, \Lambda'\Phi'=\Lambda\Phi.$
		\item $ \Lambda'\Phi^{\prime q}= \Lambda\Phi^{q}  \, \mathrm{and}\, \Lambda'\Phi'=\Lambda^q\Phi^q.$
	\end{enumerate}
	%   	\item $ \Lambda'\Phi^{\prime q}= \Lambda^q\Phi  \, \mathrm{and}\, \Lambda'\Phi'=\Lambda\Phi.$
	%   	\item $ \Lambda'\Phi^{\prime q}= \Lambda^q\Phi  \, \mathrm{and}\, \Lambda'\Phi'=\Lambda^q\Phi^q.$
	%   	\item $ \Lambda'\Phi^{\prime q}= \Lambda\Phi  \, \mathrm{and}\, \Lambda'\Phi'=\Lambda\Phi^{q}.$
	%   	\item $ \Lambda'\Phi^{\prime q}= \Lambda\Phi  \, \mathrm{and}\, \Lambda'\Phi'=\Lambda^q\Phi.$
	%   	\item $ \Lambda'\Phi^{\prime q}= \Lambda^q\Phi^q  \, \mathrm{and}\, \Lambda'\Phi'=\Lambda\Phi^{q}.$
	%   	\item $ \Lambda'\Phi^{\prime q}= \Lambda^q\Phi^q  \, \mathrm{and}\, \Lambda'\Phi'=\Lambda^q\Phi.$
	%   	\end{enumerate}
	%   
	
	For the first case,  we have  $\Lambda'\Phi^{\prime q}\Phi^{-q} =\Lambda = \Lambda'\Phi'\Phi^{-1}.$ Therefore $(\Phi'\Phi^{-1})^q=\Phi'\Phi^{-1}.$ Hence $\Phi'\Phi^{-1}$ is decomposable. Let $\Phi'\Phi^{-1}=\alpha\comp \det.$ Note that $\Lambda^{-1}\Lambda'=\Phi^{\prime-1}\Phi=(\Phi'\Phi^{-1})^{-1}=\alpha^{-1}\comp \det.$ Therefore 
	%by Theorem~\ref{thm:decomp. full}(4), we have 
	$$ \chi_{\Lambda'}^{-1} =   \chi_{(\alpha^{-1}\comp\det)\Lambda}^{-1} =  \chi_{\alpha^{-1}}^1\otimes \chi_{\Lambda}^{-1}\, \mathrm{and}\, \chi_{\Phi'}^{-1}= \chi_{(\alpha\comp \det)\Phi}^{-1} = \chi_\alpha^1\otimes \chi_{\Phi}^{-1}. $$

	For the second case,   we have $ \Lambda'\Phi^{\prime q}\Lambda^{-1}= \Phi^{q}=\Lambda'\Phi'\Lambda^{-q},$ which implies 
	$\Phi^{\prime q}\Lambda^{-1}= \Phi'\Lambda^{-q}=(\Phi^{\prime q}\Lambda^{-1})^q.$ Therefore  $\Phi^{\prime q}\Lambda^{-1}$ is decomposable. Let $\Phi^{\prime q}\Lambda^{-1}=\beta\comp\det.$ Note that $\Lambda^{\prime}\Phi^{-q} = \Phi^{\prime -q}\Lambda=(\Phi^{\prime q}\Lambda^{-1})^{-1}=\beta^{-1}\comp\det.$
	Therefore 
	%by Theorem~\ref{thm:decomp. full}(4), we have 
	$$ \chi_{\Lambda'}^{-1} =  \chi_{(\beta^{-1}\comp\det)\Phi^q}^{-1} =  \chi_{\beta^{-1}}^1\otimes \chi_{\Phi^q}^{-1}= \chi_{\beta^{-1}}^1\otimes \chi_{\Phi}^{-1}\, \mathrm{and}\, \chi_{\Phi'}^{-1}= \chi_{\Phi^{\prime q}}^{-1}= \chi_{(\beta\comp \det)\Lambda}^{-1} = \chi_\beta^1\otimes \chi_{\Lambda}^{-1}. $$

	By similar arguments, for any other case,  we can show that there exists $\gamma\in \widehat{\mathbb{F}_q^\times}$ such that either of the following holds. 
	\begin{itemize}
		\item $ \chi_{\Lambda'}^{-1} =   \chi_\gamma^1\otimes \chi_{\Lambda}^{-1}\, \mathrm{and}\, \chi_{\Phi'}^{-1}= \chi_{\gamma^{-1}}^1\otimes \chi_{\Phi}^{-1}.$
		\item $\chi_{\Lambda'}^{-1} =    \chi_\gamma^1\otimes \chi_{\Phi}^{-1}\, \mathrm{and}\, \chi_{\Phi'}^{-1}= \chi_{\gamma^{-1}}^1\otimes \chi_{\Lambda}^{-1} .$
	\end{itemize}
\end{proof}

\begin{proof}[\bf Proof of Theorem~\ref{thm:unique decomp q-1 and q+1}]
	If $\chi_{(\alpha,\beta)}^{q+1}\otimes\chi_{\Lambda}^{-1}=\chi_{(\gamma,\delta)}^{q+1}\otimes\chi_{\Phi}^{-1},$ then by restricting both sides to $Z,$ we obtain  $\alpha \beta\bar{\Lambda}=\gamma \delta\bar{\Phi}.$   
	%
	%By equating the set of  $(q+1)$-dimensional irreducible constituents   in both sides (from Theorem~\ref{thm:decomp. full}(9)), we obtain that $S_{\alpha \beta\bar{\Lambda}}^\star =S_{\gamma \delta\bar{\Phi}}^\star.$ Therefore $\alpha \beta\bar{\Lambda}=\gamma \delta\bar{\Phi}.$ 
	%Then by restricting both sides to $Z,$ we obtain  $\alpha \beta\bar{\Lambda}=\gamma \delta\bar{\Phi}.$ 
	To show the converse, suppose $\alpha \beta\bar{\Lambda}=\gamma \delta\bar{\Phi}.$ Then, we have $ \mathbb{S}_{\alpha \beta\bar{\Lambda}} =\mathbb{S}_{\gamma \delta\bar{\Phi}},$   $\mathbb{W}_{\alpha \beta\bar{\Lambda}}=\mathbb{W}_{\gamma \delta\bar{\Phi}} $ and $\{\mu \in \widehat{\mathbb{F}_q^\times} \mid \mu^2= \alpha \beta\bar{\Lambda}\}=\{\mu \in \widehat{\mathbb{F}_q^\times} \mid \mu^2= \gamma \delta\bar{\Phi}\}.$ Therefore by Theorem~\ref{thm:decomp. full}(9), we must have $\chi_{(\alpha,\beta)}^{q+1}\otimes\chi_{\Lambda}^{-1}=\chi_{(\gamma,\delta)}^{q+1}\otimes\chi_{\Phi}^{-1}.$ 
	
	Let $q\neq 3.$ To show the tensor product $\chi_{(\alpha,\beta)}^{q+1}\otimes\chi_{\Lambda}^{-1}$
	does not have unique decomposition property, let $\theta \in\widehat{\mathbb{F}_q^\times} \setminus \{1_{\mathbb{F}_q^\times}, \beta\alpha^{-1},(\beta\alpha^{-1})^2\}.$ Note that such $\theta$ exists because $q\neq 3.$ Choose a indecomposable character $\Psi \in \widehat{T_{-1}}$ such that  $\bar{\Psi}=\theta^{-1}\bar{\Lambda}.$ Then we have %$\theta\alpha\neq \beta$ and
	$(\theta\alpha)\beta\bar{\Psi}=(\theta\alpha)\beta(\theta^{-1}\bar{\Lambda})=\alpha \beta\bar{\Lambda}.$ Therefore 
	$$\chi_{(\alpha,\beta)}^{q+1}\otimes\chi_{\Lambda}^{-1}=\chi_{(\theta\alpha,\beta)}^{q+1}\otimes\chi_{\Psi}^{-1}.$$
	If $\chi_{(\theta\alpha,\beta)}^{q+1}= \chi_\gamma^1\otimes\chi_{(\alpha,\beta)}^{q+1}$ for some $\gamma \in\widehat{\mathbb{F}_q^\times},$ then either $(\theta\alpha,\beta)=(\gamma\alpha,\gamma\beta)$ or $(\theta\alpha,\beta)=(\gamma\beta,\gamma\alpha).$ Which implies either $\theta=1_{\mathbb{F}_q^\times}$ or $\theta=(\beta\alpha^{-1})^2. $ It is a contradiction to the choice of $\theta.$ Therefore 
	$\chi_{(\theta\alpha,\beta)}^{q+1}\neq  \chi_\gamma^1\otimes\chi_{(\alpha,\beta)}^{q+1}$ for all $\gamma \in\widehat{\mathbb{F}_q^\times}.$
	Since $\dim(\chi_{(\theta\alpha,\beta)}^{q+1})\neq \dim(\chi_{\Psi}^{-1}),$ we  also have $\chi_{(\theta\alpha,\beta)}^{q+1}\neq \chi_\gamma^1\otimes\chi_{\Psi}^{-1},$ for all $\gamma \in\widehat{\mathbb{F}_q^\times}.$ Hence the claim follows.
\end{proof}
\begin{remark}
L. Aburto Hageman and J. Pantoja \cite{MR1757476} also gave the description of tensor products of irreducible representation of the group $\SL_2(\mathbb{F}_q)$ in terms of induction from tori of $\SL_2(\mathbb{F}_q).$ Similar techniques can be used for complete answers of the tensor product problem of $\SL_2(\mathbb{F}_q).$
\end{remark}

\bibliography{refs}{}
\bibliographystyle{siam}
\end{document}